\documentclass[10 pt]{article}

\usepackage{enumerate}
\usepackage{fontenc,textcomp,dsfont,latexsym}
\usepackage{amsfonts,amsmath,amsthm,amssymb}
\usepackage[english]{babel}
\usepackage[latin1]{inputenc}
\usepackage{hyperref}

\def \R{\mathbb R}
\def \N{\mathbb N}

\newcommand{\F}{\mathcal{F}}

\usepackage{amssymb,amsmath,amsfonts,amsthm}
\newtheorem{Theorem}{Theorem}[part]
\newtheorem{Definition}{Definition}[part]
\newtheorem{Proposition}{Proposition}[part]

\newtheorem{Lemma}{Lemma}[part]

\newtheorem{Remark}{Remark}[part]

\newcommand{\esssup}{{\rm{ess}\;}\displaystyle{\sup}}
\newcommand{\essinf }{\displaystyle{\rm{ess\;inf}}}

\def \Ot{\; { 0\leq t\leq T}}

\def \TSO{\; \theta,S\in T_0}
\def \TS{\; \theta\in T_S}

\def \SO{\; S\in T_0}
\def \TO{\; \theta\in T_0}

\def \N{\mathbb{N}}
\def \R{\mathbb{R}}

\def\F{{\cal F}}

\def\Dzw1#1{\frac{\partial^2 #1}{\partial z \partial w_1}}

\def\Dzb1#1{\frac{\partial^2 #1}{\partial z \partial b_1}}

\begin{document}

\begin{center}
{\Large Optimal double stopping}\\
\vspace{0,5cm}
{\large Magdalena Kobylanski, Marie-Claire Quenez,  Elisabeth
Rouy-Mironescu\\}
\vspace{0,5cm}
\end{center}


\begin{abstract}
We consider the optimal double stopping time problem  defined  for each stopping time $S$ by
$v(S)=\esssup\{\, E[\psi(\tau_1, \tau_2)\, |\,\F_S], \tau_1, \tau_2  \geq S\,\}$.
Following the optimal one stopping time problem, we study the existence of optimal stopping times and give a method to compute them.
The key point is the construction of a {\em new reward} $\phi$ such that the value function $v(S)$ satisfies
$v(S)=\esssup\{\, E[\phi(\tau)\, |\,\F_S], \tau \geq S\,\}$. Finally, we give an example of an american option with double exercise time.
 
\end{abstract}

\vspace{13mm}

\noindent{\bf Key words~:}  Optimal stopping; optimal multiple stopping; stopping times; hitting times; american options with double exercise time.

\vspace{13mm}

\noindent{\bf AMS 1991 subject classifications~:} 60G40.

\section*{Introduction}

Our present work on the optimal double stopping times problem consists, following the optimal one stopping time problem, in proving the existence of the maximal reward, finding necessary or sufficent conditions for the existence of optimal stopping times, and giving a method to compute these optimal stopping times.

The results are well known in the case of the optimal one stopping time problem. Consider a 
{\em reward} given by a RCLL postive adapted process $(\phi_t, \Ot)$ on $\mathbb{F}=(\Omega,
\F,(\F_t,)_{\Ot},P)$, $\mathbb{F}$ satisfying the usual conditions, and look for the maximal 
reward 
$$v(0)= \sup \{\, E[\phi_\tau],\; \tau \in T_0\,\}\,,$$
where $T_0$ is the set of stopping times lesser than $T$.
 In order to compute $v(0)$ we introduce for each $S\in  T_0$ the {\em value function} 
 $v(S)=\esssup\{\, E[\phi_\tau\, |\,\F_S], \tau\in T_S\,\}$, where $T_S$ is the set of 
 stopping times in $T_0$ greater than $S$.
 The family $\{\,v(S), S\in T_0\,\}$ can be {\em aggregated} in a RCLL adapted process $(v_t, \Ot)$ such that $v_S=v(S)$. The process $(v_t, \Ot)$ is the {\em Snell enveloppe} of $(\phi_t, \Ot)$, that is the smallest supermartingale that dominates $\phi$.
 
 Moreover, when the reward $(\phi_t, \Ot\,\}$ is continuous an optimal stopping time is given by
 \begin{equation}\label{thetaopt}
 \theta^*(S)=\inf\{\,t\geq S, \;v_t=\phi_t\,\}\,.
 \end{equation}
 
 We show in the present work that computing the value function for the optimal double stopping times problem
 $$v(S)=\esssup\{\,E[\psi( \tau_1,\tau_2)\, |\,\F_S],\; \tau_1,\tau_2\in T_S\,\}\,,$$
 reduces to computing the value function for an optimal one stopping time problem
 $$ u(S)= \esssup\{\,E[\phi(\theta)\, |\,\F_S], \;\theta \in T_S \,\}\,,$$
 where the {\em new reward} $\phi$ is no longer a RCLL process but a family $\{\,\phi(\theta), \theta\in T_0\,\}$ of positive random variables which satisfy some compatibility properties.  
 
\vspace{0,5cm} 
 
 In  section 1, we revisit the optimal one-stopping time problem for admissible families.  In section 2, we solve the optimal two-stopping times problem. In section 3, we give an example of american exchange option with double exercise time.

\vspace{0,5cm}

Let ${\mathbb F}=(\Omega, \F, (\F_t)_{\Ot},P)$ be a probability space equipped with a filtration 
$(\F_t)_{\Ot}$ satisfying the usual conditions of right continuity and 
augmentation by the null sets of ${\cal F}= {\cal F}_T$. 
We suppose that ${\cal F}_0$ contains only sets of probability $0$ or $1$. The time horizon is a fixed constant  $T\in [0,\infty[ $. We denote by $T_{0}$ the collection of stopping times of
${\mathbb{F} }$ with values in $[0 , T]$.  More generally, for any
stopping times $S$, we denote by $T_{S}$  (resp.  $T^S$ ) the class of stopping times
$\theta\in T_0$ with $S \leq \theta$ a.s\, (resp. $\theta\leq S$ a.s). Also, we will use the following notation: we note $t_n \uparrow t$ if $\lim_{n \infty} t_n = t$ with $t_n \leq t$ for each $n$.

\section{The optimal one stopping time problem revisited}

\begin{Definition}
We say that a family 
$\{\,\phi(\theta), \TO\}$  is {\em admissible} if it satisfies  the following conditions 

\hspace{1cm}\begin{tabular}{ll}
1)& for all
$\theta\in T_0$ $\phi(\theta)$ is a $\F_\theta$-measurable positive random variable (r.v.),\\
 2)& for all
$\theta,\theta'\in T_0$, $\phi(\theta)=\phi(\theta')$ on
$\{\,\theta=\theta'\,\}$, \\
\end{tabular}
\end{Definition}

\begin{Remark}  Note that if $(\phi_t)$ is a progressive process, then the family defined by 
$\{ \phi (\theta) = \phi_{\theta}), \TO\}$ is admissible.
\end{Remark}
\noindent The {\em value function at time $S$}, where $\SO$, is given by
\begin{equation}\label{vs}
v(S)=\esssup_{\TS} E[\phi(\theta) \, |\,\F_S]\, .
\end{equation}



\begin{Proposition}\label{P1.SuperM}
Let $\{\,\phi(\theta), \;\theta\in T_0\,\}$ be an admissible family, 
then the family of r.v $\{\, v(S) , \; S \in T_0 \,\}$ is admissible and is a {\em supermartingale} system, that is 
for any stopping times 
$\theta, \theta^{'}$ $ \in$ $T_S$ such that $\theta \geq \theta^{'}$, 
$E[v (\theta) \, |\,{\cal F}_{  \theta^{'} }]  \leq  v (\theta^{'}) \quad
\mbox{a.s\,,} $.

\end{Proposition}




\begin{Definition}
An admissible family $\{\,\phi(\theta),\TO\,\}$ is said to be  {\em right (resp. left) continuous along stopping times in expectation } if for any  $\theta\in T_0$ and for any $\theta_n \downarrow \theta$  (resp. 
$\theta_n \uparrow \theta$) a.s. one has
$ \displaystyle{E[\phi(\theta)]=\lim_{n \to\infty} E[\phi(\theta_n)].} $
\end{Definition}

Recall the following classical lemma (see El Karoui (1979))
\begin{Lemma}\label{L1}
Let $\{\,\phi(\theta),\TO\,\}$ be an admissible family right continuous along stopping times in expectation 
such that $E [\esssup_{\theta\in
T_0}\phi(\theta)]<\infty$.
Then, the family $\{\, v(S) , \; S \in T_0 \,\}$ is right continuous along stopping times in expectation.
\end{Lemma}

We will now state the existence of of an optimal stopping time under quite weak assumptions.\\
For each  $S$ $\in$ $T_0$, let us introduce the following $\F_S$-measurable random variable 
 $\theta^*(S)$ defined by 
\begin{equation}\label{thetaop}\theta^*(S) :=\essinf \{\,\theta \in T_S\, , \, v(\theta) =\phi(\theta) \,\,\mbox {a.s} \,\}
\end{equation}

Note that $\theta^*(S)$ is a stopping time. Indeed, for each $S$ $\in$ $T_0$, 
one can easily show that the set ${\cal T}_S= \{\,\theta \in T_S\, , \, v(\theta) =\phi(\theta) \,\,\mbox {a.s} \,\}$ is closed under pairwise maximization. By a classical result, there
exists a sequence $(\theta^n) _{n \in \N}$ of stopping times in ${\cal T}_S$ such that 
$\theta_n \downarrow \theta^*(S)$ a.s. Futhermore, we state the following theorem which generalizes  the classical existence result of optimal stopping times to the case of a reward given by an admissible 
family of r.v. (instead of a RCLL adapted process).

\begin{Theorem}\label{T1.2}
Suppose that $E [\esssup_{\theta\in
T_0}\phi(\theta)]<\infty$.

Let $\{\, \phi(\theta), \theta\in T_0\,\}$ be an admissible family, right and left  continuous along 
stopping times in expectation. Let $\{\,v(S), \SO\,\}$ be the family of value function defined by (\ref{vs}).
Then, for each $S$ $\in$ $T_0$, the stopping time $\theta^*(S)$ defined by (\ref{thetaop}) 
is optimal for $v(S)$ {\em (i.e}\quad $v(S)= E[\phi(\theta^*(S))\, |\,\F_S]\;{\rm )}$.
\end{Theorem}

\begin{Remark} When $\phi$ is given by a (RCLL) right continuous left limited adapted process, since the value function can be aggregated by a RCLL process $(v_t)_{t \in[0,T]}$, $\theta^*(S)$ satisfies equality (\ref{thetaopt}) a.s.
\end{Remark}
\noindent{\sc Sketch of the Proof}~:  Fix $S$ $\in$ $T_0$. 
As in the case of a reward process, we begin by constructing a family of stopping times that are approximatively optimal. For 
$\lambda$ $\in$ $]0,1[$, define the $\F_S$-measurable random variable 
$\theta^{\lambda} (S)$ by

\begin{equation} \label{thetal}
\theta^{\lambda} (S) := \essinf\{\,\theta \in T_S\, , \,\lambda v(\theta) \leq \phi(\theta) \,\,\mbox {a.s} \,\} 
\end{equation}

For all $\SO$, the function $\lambda \mapsto \theta^\lambda(S)$ is increasing on $\in]0,1[$ (and bounded above by $\theta^*(S)$). Let us now define, for $\SO$,
$\hat \theta(S):= \lim_{\lambda\uparrow 1} \uparrow \theta^\lambda(S)\,.$\\
Then , we  prove that $\hat \theta(S)$ is optimal for $v(S)$ and that $\hat\theta(S)=\theta^*(S)$ by similar arguments as in the case of a reward process but the proof is shorter since it is no longer necessary to aggregate the value function. 
$\Box$


\section{The optimal two stopping times problem}

We consider now the optimal two stopping times problem. 
Let us introduce the following definitions

\begin{Definition} The family 
$\{\,\psi(\theta,S),\TSO\,\}$ is  an {\em biadmissible family } if it satisfies\vspace{5pt}

 \noindent \hspace{1cm}\begin{tabular}{ll}
 1)& for all $ \theta,S\in T_0,$  $\psi(\theta,S)$ is a $\F_{\theta\vee S}$-measurable positive random variable,\\
 2)&  for all $\theta,\theta',S,S'\in
T_0,\quad \psi(\theta,S)=\psi(\theta',S')$ a.s on $
\{\,\theta=\theta'\,\}\cap\{\,S=S'\,\} $, \\
\end{tabular}
\end{Definition}

\vspace{0,5cm}

For each stopping time $S$, let us consider the value function
\begin{equation}\label{vS}
v(S)=\esssup_{\tau_1,\tau_2\in T_S} E[\psi(\tau_1,\tau_2)\, |\,\F_S]
\end{equation}

As in the case of one stopping time problem,  the family $\{\,v(S), \SO\,\}$ 
is an admissible family of postive r.v.  and the family of positive r.v $\{\, v(S) , \SO \,\}$ is a supermartingale system.

\subsection{Reduction to an optimal one stopping time problem}

This optimal two-stopping times problem can be expressed in terms of 
one-stopping time problems as follows.
 \noindent For each stopping
time $\theta \in T_S$, define the two $\F_{\theta}$-measurable
random variables
\begin{equation}\label{u12}
u_1({\theta})=\esssup_{\tau_1\in T_{\theta}} E[\psi(\tau_1,
\theta)\, |\,\F_{\theta}], \quad  u_2({\theta}) =\esssup_{\tau_2\in
T_{\theta}} E[\psi(\theta,\tau_2)\, |\,\F_{\theta}].\quad
\end{equation}
Note that by Proposition \ref{P1.SuperM} since $\{\,\psi(\theta,S),\TSO\,\} $ is admissible, the 
families $\{\,u_1(\theta), \TS\,\}$ and $\{\,u_2(\theta),\TS\,\}$ are admissible.
Put
\begin{equation}\label{phi}
\phi(\theta)=\max[u_1(\theta), u_2(\theta)],
\end{equation}
the family $\{\,\phi(\theta), \TS\,\}$, which is called the {\em new reward family}, is also clearly admissible.
\noindent
Define the $\F_S$-measurable variable
\begin{equation}\label{uS}
u(S)=\esssup_{\theta\in T_S} E[\phi(\theta)\, |\,\F_S]\quad
\mbox{a.s}\; .
\end{equation}

\begin{Theorem} \label{T3}{\em (Reduction)}
Suppose that $\{\,\psi(\theta,S), \TSO\,\}$ is an admissible family. For each stopping time $S$, consider $v(S)$ defined
by (\ref{vS}) and $u(S)$ defined by (\ref{u12}), (\ref{phi}),
(\ref{uS}), then $v(S)=u(S) \quad \mbox {a.s\,.}$
\end{Theorem}
\noindent{\sc Proof} ~: Let $S$ be a stopping time.\\
\underline{Step 1~:} Let us show first that
$
 v(S)\leq u(S) \quad \mbox{a.s}\, .
$\\
 Let $\tau_1,\tau_2\in
T_S$. Put $A=\{\,\tau_1\leq \tau_2\,\}$. As $A$ is in $\F_{\tau_1\wedge \tau_2}$,
$$E[\psi(\tau_1,\tau_2)\, |\,\F_S]= E[
{\bf 1}_AE[\psi(\tau_1,\tau_2)\, |\,\F_{\tau_1\wedge \tau_2}]\, |\,\F_S]
+E[{\bf 1}_{A^ c}E[\psi(\tau_1,\tau_2)\, |\,\F_{\tau_1\wedge \tau_2}]\, |\,\F_S].
$$

\noindent By noticing that on  $A$ one has
$
E[\psi(\tau_1,\tau_2)\, |\,\F_{\tau_1\wedge
\tau_2}]=E[\psi(\tau_1,\tau_2)\, |\,\F_{\tau_1}] \leq
u_2({\tau_1})\leq \phi({\tau_1\wedge \tau_2}) \quad \mbox{a.s},
$
and similarly on $A^ c$ one has
$E[\psi(\tau_1,\tau_2)\, |\,\F_{\tau_1\wedge
\tau_2}]=E[\psi(\tau_1,\tau_2)\, |\,\F_{\tau_2}] \leq
u_1({\tau_2})\leq \phi({\tau_1\wedge \tau_2})\quad \mbox{a.s},
$
 we get
$
E[\psi(\tau_1,\tau_2)\, |\,\F_S]
\leq E[ \phi({\tau_1\wedge \tau_2}) \, |\,\F_S]\leq u(S)\quad
\mbox{a.s}\; .
$\\
 By taking the supremum over $\tau_1$ and $\tau_2$ in $T_S$
we obtain step 1.

\noindent 
\underline{Step 2~:} To simplify, suppose the existence of optimal stopping times for the three one-optimal stopping times problems.
Let $\theta^{*}$ be an optimal stopping time for $u(S)$ {\em i.e}
\begin{equation*}\label{theta*}u(S)=\esssup_{\theta\in T_S}
 E[\phi({\theta})\, |\,\F_{S}] = E[\phi({\theta^{*}})\, |\,\F_{S}]\quad
 \mbox{a.s}\,,
 \end{equation*}
 and let $\theta^*_1$ and
 $\theta^*_2$ be optimal stopping times for $u_1({\theta^*})$ and $u_2({\theta^*})$,  i.e.
 $u_1({\theta^*})= E[\psi(\theta^*_1,\theta^{*})\, |\,\F_{\theta^*}]$ a.s.  and   $u_2({\theta^*})=
 E[\psi(\theta^{*}, \theta^*_2)\, |\,\F_{\theta^*}]$  a.s.\\
  Put $B=\{\,\; u_1(\theta^*) \leq u_2(\theta^ *) \; \,\}$. Let $(\tau_1^*,\tau_2^*)$ be the stopping times defined by 
\begin{equation}\label{tau12}
 \tau_1^*=\theta^{*}{\bf 1 }_{B}+ \theta^{*}_1{\bf 1 }_{B^c}; \quad 
 \tau_2^*=\theta^*_2{\bf 1 }_{B}+ \theta^{*}{\bf 1 }_{B^c}, 
 \end{equation}

\noindent As $B$
is $\F_{\theta^*}$ measurable, and
 since $u({S}) =
E[\phi({\theta^{*}})\, |\,\F_{S}]$, we have
\begin{eqnarray*}
u({S})& =& E[{\bf 1}_B  u_2({\theta^*}) + {\bf 1 }_{B^c}
u_1({\theta^*})\, |\,\F_S].
\end{eqnarray*}
Hence, we have
$u({S}) = E[{\bf 1}_B \psi(\theta^{*} ,\theta^*_2) + {\bf 1
}_{B^c}\psi(\theta^*_1 ,\theta^{*})\, |\,\F_S] =E[\psi(\tau_1^*,\tau_2^*)\, |\,\F_S]
$.

Hence,  $u({S})\leq v(S)$ a.s.  Note that in the general case, Step 2 can be also be obtained, for example by using optimizing sequences of stopping times.
$\quad$ $\Box$
 
\begin{Remark} \label{remarque}
Note that step 2 gives a constructive method of optimal stopping times. Indeed, let 
$\theta^ *$ be optimal for $u(S)$, let
$\theta^*_2$ be optimal for $u_2({\theta^*})$ and let $\theta^*_1$ be optimal for
$u_1({\theta^*})$.
Put $B=\{\,\; u_1(\theta^*) \leq u_2(\theta^ *) \; \,\}$, then the pair of stopping times $(\tau_1^*,\tau_2^*)$ defined by (\ref{tau12}) is optimal for $v(S)$. Note that step 1 gives also a necessary condition of optimality (but not a characterization in general).
\end{Remark}

Before studing  the problem of existence of optimal stopping times, we have to state some regularity properties of the new reward family $\{\phi(\theta),\TO\}$. 
\subsection{Regularity properties of the new reward family}
As in the case of one stopping time, we will need some regularity properties on the reward family $\psi$.
Let us introduce the following definition,
\begin{Definition} 
Suppose that $E [\esssup_{\theta, S \in
T_0}\psi (\theta, S)]<\infty$.\\
An admissible family $\{\,\psi(\theta,S),\TSO\,\}$ is said to be {\em uniformly right (resp. left ) continuous in expectation along stopping times } if 
for each $\theta$, $S$ $\in$ $T_0$ and each sequences $(\theta_n)$, $(S_n)$ of $T_0$ such that 
$\theta_n \downarrow \theta$ a.s. and $S_n \uparrow S$ a.s., we have 
$\lim_{n \infty} \,\,E\left( \esssup_{\theta \in T_0 }   \vert \psi (\theta, S) - \psi (\theta , S_n) \vert \right) = 0$
 and  
$\lim_{n \infty} \,\, E\left( \esssup_{\theta \in T_0  }   \vert \psi (S, \theta) - \psi (S_n , \theta) \vert \right)=0.
$ ( resp. 
$\lim_{n \infty} \,\,E\left( \esssup_{S \in T_0 }   \vert \psi (\theta, S) - \psi (\theta_n , S) \vert \right)$ $=$ $0$
 and
$\lim_{n \infty} \,\, E\left( \esssup_{S \in T_0  }   \vert \psi (S, \theta) - \psi (S , \theta_n) \vert \right)$ $=$ $0$).

\end{Definition}

The following continuity property holds true for the new reward family~:
\begin{Theorem} \label{T2.4}
Suppose that the admissible family $\{\,\psi(\theta,S),\TSO\,\}$ is uniformly right continuous(resp. left continuous) in expectation along stopping times. \\
 Then, the  family $\{\,\phi(S), \SO\,\}$  defined by (\ref{phi}) is right continuous (resp. left continuous) in expectation along stopping times.

\end{Theorem}
\noindent{\sc Proof~:} Let us show the right continuity property.
As $\phi(\theta)=\max[u_1(\theta),u_2(\theta)]$, it is sufficent to show the right continuity property for the family $\{\,u_1(\theta),\TO\,\}$.\\
Let us introduce the following value function for each $S,\theta$ $\in$ $T_0$,
\begin{equation}\label{U1}
U_1 (\theta, S) = \esssup_{\tau_1 \in T_{\theta}} E[\psi (\tau_1, S)\, |\,\F_{\theta}]\quad
\mbox{a.s}\,.
\end{equation}
As for all $\TO$, $u_1(\theta)=U_1(\theta,\theta)\quad \mbox{a.s}\,,$
it is sufficient to prove that $\{\,U_1 (\theta, S),   \theta, S \in T_0\,\}$ is right continuous in expectation along stopping times.\\
Let $\theta\in T_0$ and $(\theta_n)_n$ a nonincreasing sequence of stopping times in $T_{\theta}$ that converges to $\theta$ a.s\,. 
Let $S \in T_0$ and $(S_n)_n$ a nonincreasing sequence of stopping times in $T_S$ that converges to $S$ a.s\,. 

We have,
$$
|  E[U_1 (\theta, S)] - E[U_1 (\theta_n, S_n)]  | \quad\leq \quad | E[U_1 (\theta, S)] - E[U_1 (\theta_n , S)] | + | E[U_1 (\theta_n, S)] - E[U_1 (\theta_n, S_n)] .
$$

The first term of the right hand side converges a.s  to $0$ as $n$ tends to $\infty$ by the right continuity in expectation of  the value function $\{U_1 (\theta, S), \theta \in T_0\}$ by lemma \ref{L1}.

Also, the seond term converges to $0$ since  $\{\,\psi(\theta,S),\TSO\,\}$ is uniformly right continuous in expectation along stopping times.

For the left continuity property of the new reward, the arguments are the same. \quad $\Box$



We will now turn to the problem of existence of optimal stopping times.
\subsection{Existence of optimal stopping times}
 Under uniform right and left  continuity in expectation along stopping times of the admissible family $\{\,\psi(\theta,S),$ $\TSO\,\}$, we can establish the existence of optimal
stopping times. 

\begin{Theorem}\label{Topt} 
Let $\{\,\psi(\theta,S), \TSO\,\}$ be an admissible family which is uniformly right and left  continuous in expectation along stopping times.
Suppose that there exists $p>1$ such that
 $\displaystyle{
E[\esssup_{\theta,S\in T_0}\psi(\theta,S)^p]<\infty}$.\\
 Then, there exists a pair $(\tau_1^*,\tau_2^*)$ of optimal stopping times for $v(S)$ defined by (\ref{vS}), that is
$$v(S)= \esssup_{\tau_1,\tau_2\in T_S} E[\psi(\tau_1,\tau_2)\,|\F_S\,]= E[\psi(\tau_1^*,\tau_2^*)\,|\F_S\,]\,.$$

\end{Theorem}
\noindent{\sc Proof ~:} By Theorem \ref{T2.4}  the admissible family of positive r.v $\{\,\phi(\theta), \TO\,\}$ defined by (\ref{phi}) is right and left  continuous in expectation along stopping times. By Theorem \ref{T1.2} the stopping time
\[\theta^{*}(S)=\essinf\{\,\theta \in T_S\; | \; v(\theta)=\phi(\theta)\, \mbox{a.s} \,\},\]
is optimal for $u(S)= v(S)$,
that is
$u(S)=\esssup_{\theta\in T_S} E[\phi({\theta})\, |\,\F_{S}]  =
E[\phi(\theta^{*})\, |\,\F_{S}]\quad \mbox{a.s}\;.$
Moreover, the families  $\{\,\psi(\theta,\theta^*),\; \theta\in T_{\theta^*}\,\}$ and $\{\,\psi(\theta^*,\theta),\; \theta\in T_{\theta^*}\,\}$ 
are clearly admissible and are right and left  continuous in expectation along stopping times. Consider the following optimal stopping time problem~:
$$ v_1(S)=\esssup_{\TS}E[\psi(\theta,\theta^*)\,|\, \F_{S}] \quad \mbox{and}\quad 
v_2(S)=\esssup_{\TS}E[\psi(\theta^*,\theta)\,|\, \F_{S}].$$
 Let  $\theta^*_1$ and $\theta^*_2$ 
are optimal stopping times given by Theorem \ref{T1.2} for $v_1(\theta^*)$ and
$v_2(\theta^*)$. Note that  $v_1(\theta^*)=u_1(\theta^*)$ and
$v_2(\theta^*)=u_2(\theta^*)$.\\
The existence of optimal stopping times
$(\tau^*_1,\tau^*_2)$ for problem (\ref{vS}) follows  by Remark \ref{remarque}. \quad $\Box$

\section{Example:  an exchange option}
Suppose that the market contains two risky assets with price processes
$X_t= (X^1_t, X^2_t)$ which satisfy 
$dX^1_t = rX^1_t dt + \sigma_1 X^1_t dW^1_t$ and 
$dX^2_t = rX^2_t dt + \sigma_2 X^2_t dW^2_t,$
where $W= (W^1,W^2)$ is a $\R^2$-valued brownian motion , and
$(\F_t)_{t\geq 0}$ is the associated standard filtration. Without loss of generality, we first suppose that the interest rate $r$ is equal to $0$ (the general case can be derived by 
considering discounted prices).
The buyer of the option chooses two
stopping times $\tau_1$ (for the first asset) and $\tau_2$ (for the
second one) smaller than $T$. At time $\tau_1
\vee \tau_2$, he is allowed to exchange the first asset against the second one; in
other terms, he receives the amount:
$(X^1_{\tau_1}- X^2_{\tau_2})^+.$

Given the initial data $(s,x)= (s, x_1, x_2) \in[0,T]\times \R^2$, the price
at time $s$ of this option is given by
\begin{eqnarray} \label{vsx}
v(s,x)={\rm ess} \sup_{ \hspace{-0,5cm} \tau_1,\tau_2\in T_s} E[(X^{1,s,x_1}_{\tau_1}- X^{2,s,x_2}_{\tau_2})^+|\F_s],\;\; &&
\end{eqnarray}

Given initial data $(s,x_1)\in[0,T]\times \R$ consider also
${\displaystyle v^1(s,x_1,y_2) ={\rm ess} \sup_{ \hspace{-0,5cm} \tau_1\in T_s}  E _{s,x_1}[
(X^{1}_{\tau_1}- y_2)^+]},$ 
for each fixed parameter $y_2 \in  \R$. We have clearly that
$ v^1(s,x_1,y_2) =y_2 C_1 (s, \frac{x_1}{y_2}),$
where for initial data $(s,z)\in[0,T]\times \R$,
$ C_1 (s, z) = E _{s,z} [(X^{1}_T- 1) ^+],$ which is explicitely given by
the Black-Sholes formula.
 Similarly, given initial data $(t,y_2)\in[0,T]\times \R$,
consider \\
${\displaystyle v^2(x_1,t,y_2) ={\rm ess} \sup_{ \hspace{-0,5cm} \tau_2\in T_t}  E _{t,y_2}[(x_1 - X^{2}_{\tau_2})^+]}$
for each fixed parameter $x_1 \in \R$.
 We have clearly that
$
v^2(x_1,t,y_2) =x_1 P_2 (t, \frac{y_2}{x_1}),
$
where for initial data $(s,z)\in[0,T]\times \R$,
$ P_2 (t, z) =  E _{t,z}[(1- X^{2}_T) ^+].$ Note that the function $P_2 (s, z)$ corresponds to the price of an european
put option. We set
$$\phi(s,x_1, x_2)=\max{\left( x_2 C_1(s, \frac{x_1}{x_2}), x_1 P_2(s,
\frac{x_2}{x_1})\right)},
\mbox{ for }s,x_1, x_2
\in[0,T]\times\R^2.$$

By the previous results, the price function $v$ is the value function
of the following optimal stopping time problem:
\begin{eqnarray*}\label{ee}
    v(s,x_{1}, x_{2}) = \esssup_{\theta \in T_s} E_{s, x_{1}, x_{2}}[ \phi(\theta, X^{1}_{\theta}, X^{2}_{\theta})].,\;\; &&
\end{eqnarray*}
It corresponds to the price of an American option with maturity $T$ and payoff $\left( \phi(t, X^{1}_t, X^{2}_t) \right)_{0\leq t \leq T}$.
By classical results, it is a viscosity solution of an obstacle problem. By Remark \ref{remarque}, we can easily construct a pair $(\tau_1 ^*, \tau_2 ^*)$ of optimal stopping times for the exchange option 
(\ref{vsx}).

\end{document}